\documentstyle[11pt]{article}
\begin{document}
\title{{\bf Global Solution of Atmospheric Circulation Models with Humidity Effect}
\thanks{Foundation item: the National Natural Science Foundation of China
(NO. 11271271).}}
\author{{\sl
Hong Luo}
\\{\small College of Mathematics and Software Science,
Sichuan Normal University,} \\
 {\small Chengdu, Sichuan 610066, China }}
\date{}
\maketitle
 \begin{minipage}{5.5 in}
  \noindent{\bf Abstract:}\,\ \ The atmospheric circulation models are
  deduced from the very complex atmospheric circulation models
 based on the actual background and meteorological data. The models are able to show features of
  atmospheric circulation and are easy to be studied.
  It is proved that existence of global solutions to atmospheric
 circulation models
with the use of the $T$-weakly continuous operator.


\end{minipage}

\section{Introduction}
Mathematics is a summary and abstraction of the production practices
and the natural sciences, and is a powerful tool to explain natural
phenomena and reveal the laws of nature as well. Atmospheric
circulation is one of
 the main factors affecting the global climate, so it is very necessary
 to understand and master its mysteries and laws. Atmospheric circulation
 is an important mechanism to complete the transports and balance of atmospheric
 heat and moisture and the conversion between various energies. On the contrary,
 it is also the important result of these physical transports, balance and conversion.
 Thus it is of necessity to study the characteristics, formation, preservation, change
 and effects of the atmospheric circulation and master its evolution law, which is not
 only the essential part of human's understanding of nature, but also the helpful method
 of changing and improving the accuracy of weather forecasts, exploring global climate change,
 and making effective use of climate resources.
\par
The atmosphere and ocean around the earth are rotating geophysical
fluids, which are also two important components of the climate
system. The phenomena of the atmosphere and ocean are extremely rich
in their organization and complexity, and a lot of them cannot be
produced by laboratory experiments. The atmosphere or the ocean or
the coupled atmosphere and ocean can be viewed as an initial and
boundary value problem \cite{Phillips} \cite{Rossby}, or an infinite
dimensional dynamical system
\cite{Lions1}\cite{Lions2}\cite{Lions3}. These phenomena involve a
broad range of temporal and spatial scales \cite{Charney2}. For
example, according to \cite{von Neumann}, the motion of the
atmosphere can be divided into three categories depending on the
time scale of the prediction. They are motions corresponding
respectively to the short time, medium range and long term behavior
of the atmosphere. The understanding of these complicated and
scientific issues necessitate a joint effort of scientists in many
fields. Also, as \cite{von Neumann} pointed out, this difficult
problem involves a combination of modeling, mathematical theory and
scientific computing.

Some authors have studied the atmospheric motions viewed as an
infinite dimensional dynamical system. In \cite{Lions1}, the authors
study as a first step towards this long-range project which is widely
considered as the basic equations of atmospheric dynamics in
meteorology, namely the primitive equations of the atmosphere. The
mathematical formulation and attractors of the primitive equations,
with or without vertical viscosity, are studied. First of all, by
integrating the diagnostic equations they present a mathematical
setting, and obtain the existence and time analyticity of solutions
to the equations. They then establish some physically relevant
estimates for the Hausdorff and fractal dimensions of the attractors
of the problems. In \cite{Li}, based on the complete dynamical
equations of the moist atmospheric motion, the qualitative theory of
nonlinear atmosphere with dissipation and external forcing and its
applications are systematically discussed by new theories and
methods on the infinite dimensional dynamical system. In
\cite{Wang}, by Lions theorem in the Hilbert space, the existence
and uniqueness of the weak solution of water vapour-equation with
the first boundary-value problem are proven, and the scheme of the
finite-element method according to the weak solution is proposed. In
\cite{Huang}, the model of climate for weather forecast is studied,
and the existence of the weak solution is proved by Galerkin method.
The asymptotic behaviors of the weak solution is described by the
trajectory attractors.

In this article, Atmospheric circulation equations with humidity
effect is considered, which is different from the previous research.
Previous studies are based on two kinds of equations, one is about
the heat and humidity transfer(\cite{Li},\cite{Wang}), without
considering the diffusion of heat and of humidity; the other is
p-coordinates equation in(\cite{Lions1},\cite{Huang}), being used to
consider the horizontal movement of the atmosphere, and the vertical
direction of the velocity is transformed into pressure and
topography. Atmospheric circulation equations with humidity effect,
derivation from the original equations in \cite{Richardson},
considering the diffusion of heat and of humidity, not be restricted
by p-coordinates, can be deformed based on the different concerns.
In the last part of the article, Existence of global solutions to
the atmospheric circulation models is obtained, which implies that
atmospheric circulation has its own running way as humidity source
and heat source change, and confirms that the atmospheric
circulation models are reasonable.

The paper is organized as follows. In Section 2 we present
derivation of atmospheric circulation models. In Section 3, we prove
that the atmospheric circulation models possess global solutions by
using space sequence method.

\section{Derivation of Atmospheric Circulation Models with \\Humidity Effect}
\subsection{Original Model}
The hydrodynamical equations governing the atmospheric circulation
are the Navier-Stokes equations with the Coriolis force generated by
the earth¡¯s rotation, coupled with the first law of thermodynamics.
\par
Let $(\varphi, \theta, r)$ be the spheric coordinates, where
$\varphi$ represents the longitude, $\theta$ the latitude, and $r$
the radial coordinate. The unknown functions include the velocity
field $u=(u_\varphi, u_\theta, u_r)$, the temperature function $T$,
the humidity function $q$, the pressure $p$ and the density function
$\rho$. Then the equations governing the motion and states of the
atmosphere consist of the momentum equation, the continuity
equation, the first law of thermodynamic, and the diffusion equation
for humidity, and the equation of state (for ideal gas):
\begin{eqnarray}
\label{eqa1} \rho[\frac{\partial u}{\partial t}+\nabla_u
u+2\overrightarrow{\Omega} \times u]+\nabla
p+\overrightarrow{\kappa} \rho g=\mu \Delta u,
\end{eqnarray}
\begin{eqnarray}
\label{eqa2} \frac{\partial \rho}{\partial t}+div(\rho u)=0,
\end{eqnarray}
\begin{eqnarray}
\label{eqa3} \rho c_v[\frac{\partial T}{\partial t}+u\cdot\nabla T]
+p div u=Q+\kappa_T \Delta T,
\end{eqnarray}
\begin{eqnarray}
\label{eqa4} \rho[\frac{\partial q}{\partial t}+u \cdot
\nabla_q]=G+\kappa_q \Delta q,
\end{eqnarray}
\begin{eqnarray}
\label{eqa5} p=R_0\rho T,
\end{eqnarray}
Where $-\infty <\varphi <+\infty$, $-\frac{\pi}{2}<\theta
<\frac{\pi}{2}$, $r_0<r<r_0+h$, $r_0$ is the radius of the earth,
$h$ is the height of the troposphere, $\Omega$ is the earth's
rotating angular velocity, $g$ is the gravitative constant, $\mu,
\kappa_T, \kappa_q, c_v, R_0$ are constants, $Q$ and $G$ are heat
and humidity sources, and $\overrightarrow{\kappa} = (0, 0, 1)$. The
differential operators used are as follows:
\par
1. The gradient and divergence operators are given by:
$$
\nabla=(\frac{1}{r\cos\theta}\frac{\partial}{\partial \varphi},
\frac{1}{r} \frac{\partial}{\partial \theta},
\frac{\partial}{\partial r}),
$$
$$
div u=\frac{1}{r^2}\frac{\partial}{\partial r}(r^2 u_r)+\frac{1}{r
\cos\theta} \frac{\partial (u_\theta \cos \theta)}{\partial
\theta}+\frac{1}{r \cos\theta} \frac{\partial u_\varphi}{\partial
\varphi},
$$
\par
2. In the spherical geometry, although the Laplacian for a scalar is
different from the Laplacian for a vectorial function, we use the
same notation $\Delta$ for both of them:
$$
\Delta u=( \Delta u_\varphi+\frac{2}{r^2 \cos\theta} \frac{\partial
u_r }{\partial \varphi}+\frac{2 \sin \theta}{r^2 \cos^2 \theta}
\frac{\partial u_\theta}{\partial \varphi}- \frac{u_\varphi}{r^2
\cos^2 \theta},
$$
$$
\Delta u_\theta+\frac{2}{r^2} \frac{\partial u_r }{\partial
\theta}-\frac{u_\theta}{r^2 \cos^2 \theta}- \frac{2 \sin\theta}{r^2
\cos^2 \theta}\frac{u_\varphi}{\partial \varphi},
$$
$$
 \Delta u_r+\frac{2 u_r}{r^2}+\frac{2}{r^2 \cos \theta}
\frac{\partial (u_\theta \cos\theta)}{\partial \theta}- \frac{2}{r^2
\cos \theta}\frac{u_\varphi}{\partial \varphi}),
$$
$$
\Delta f=\frac{1}{r^2 \cos\theta} \frac{\partial
f}{\partial\theta}(\cos\theta \frac{\partial}{\partial
\theta})+\frac{1}{r^2 \cos^2 \theta} \frac{\partial^2 f}{\partial
\varphi^2}+\frac{1}{r^2} \frac{\partial f}{\partial r}(r^2
\frac{\partial}{\partial r}),
$$
\par
3. The convection terms are given by
$$
\nabla_u u=(u\cdot \nabla u_\varphi+\frac{u_\varphi
u_r}{r}-\frac{u_\varphi u_\theta}{r} \tan \theta,
$$
$$
u\cdot \nabla u_\theta+\frac{u_\theta u_r}{r}+\frac{u^2_\varphi}{r}
\tan \theta, u\cdot \nabla u_r-\frac{u^2_\varphi+ u^2_\theta}{r}),
$$
\par
4. The Coriolis term $2 \overrightarrow{\Omega}\times u$ is given by
$$
2 \overrightarrow{\Omega}\times u=2 \Omega(\cos \theta
u_r-\sin\theta u_\theta, \sin \theta u_{\theta}, -\cos\theta u_\varphi),
$$
Here $\overrightarrow{\Omega}$ is the angular velocity vector of the
earth, and $\Omega$ is the magnitude of the angular velocity.
\par
They are supplemented with the following initial value conditions
\begin{eqnarray}
\label{eqa100} (u, T, q ) = (\varphi_{10},\varphi_{20},
\varphi_{30}) \quad at \quad t = 0.
\end{eqnarray}
\par
 Boundary conditions are needed at the
top and bottom boundaries $(r_0, r_0+h)$. At the top and bottom
boundaries $(r=r_0,r_0+h)$, either the Dirichlet boundary condition
or the free boundary condition is given
\begin{eqnarray}
\label{eqa101} (\hbox{Dirichlet})\quad \left\{
   \begin{array}{ll}
(u, T,q ) =(0, T_0, q_0), & r=r_0,\\
(u, T,q )=(0, T_1, q_1), & r=r_0+h,
\end{array}
\right.
\end{eqnarray}
\begin{eqnarray}
\label{eqa102} (\hbox{free})\quad \left\{
   \begin{array}{ll}
(u, T,q ) =(0, T_0, q_0), & r=r_0,\\
(u_r,T,q )=(0,T_1, q_1),\quad \frac{\partial(u_\varphi,
u_\theta)}{\partial r} =
 0 & r=r_0+h,
\end{array}
\right.
\end{eqnarray}
\par
For $\varphi$, periodic condition are usually used,
 for any integers $k_1\in Z$
\begin{eqnarray}
\label{eqa103} (u, T, q )(\varphi + 2k_1 \pi, \theta , r) = (u, T, q
)(\varphi + 2k_1 \pi, \theta , r).
\end{eqnarray}
\par
Because $(\varphi, \theta, r)$ are in a circular field with
$C^\infty$ boundary, the space domain is taken as $(0,2\pi) \times
(-\frac{\pi}{2}, \frac{\pi}{2}) \times (r_0,r_0+h)$ and periodic
condition is written as
\begin{eqnarray}
\label{eqa104} (u, T, q )(0,\theta, r) = (u, T, q )(2\pi,\theta, r).
\end{eqnarray}

 For simplicity, we
study the  problem with the Dirichlet boundary conditions, and all
results hold true as well for other combinations of boundary
conditions. Atmospheric convection equations can be read as
(\ref{eqa1})-(\ref{eqa101}), (\ref{eqa104}).

\par
The above equations were basically the equations used by L. F.
Richardson in his pioneering work \cite{Richardson}. However, they
are in general too complicated to conduct theoretical analysis. As
practiced by the earlier workers such as J. Charney, and from the
lessons learned by the failure of Richardson's pioneering work, one
tries to be satisfied with simplified models approximating the
actual motions to a greater or lesser degree instead of attempting
to deal with the atmosphere in all its complexity. By starting with
models incorporating only what are thought to be the most important
of atmospheric influences, and by gradually bringing in others, one
is able to proceed inductively and thereby to avoid the pitfalls
inevitably encountered when a great many poorly understood factors
are introduced all at once. The simplifications are usually done by
taking into consideration of some main characterizations of the
large-scale atmosphere. One such characterization is the small
aspect ratio between the vertical and horizontal scales, leading to
hydrostatic equation replacing the vertical momentum equation. The
resulting system of equation are called the primitive equations
(PEs); see among others \cite{Lions1}. The another characterization
of the large scale motion is the fast rotation of the earth, leading
to the celebrated quasi-geostrophic equations \cite{Charney1}.
\par
 For convenience of research, the approximations we adopt
involves the following components:
\par
First, we often use Boussinesq assumption, where the density is
treated as a constant except in the buoyancy term and in the
equation of state.
\par
Second, because the air is generally compressible, we do not
use the equation of state for ideal gas, rather, we use the
following empirical formula, which can by considered as the linear
approximation of
\begin{eqnarray}
\label{eqa6} \rho=\rho_0[1-\alpha_T(T-T_0)+\alpha_q(q-q_0)],
\end{eqnarray}
where $\rho_0$ is the density at $T = T_0$ and $q = q_0$, and
$\alpha_T$ and $\alpha_q$ are the coefficients of thermal and
humidity expansion.
\par
Third, since the aspect ratio between the vertical scale and the
horizontal scale is small, the spheric shell the air occupies is
treated as a product space $S^2_{r_0}\times(r_0, r_0+ h)$. This
approximation is extensively adopted in geophysical fluid dynamics.
Under the above simplification, we have the following equations
governing the motion and states of large scale atmospheric
circulations:
\begin{eqnarray}
\label{eqa7} \frac{\partial u}{\partial t}+\nabla_u u=\nu \Delta
u-2\overrightarrow{\Omega} \times u-\frac{1}{\rho_0}\nabla
p-[1-\alpha_T(T-T_0)+\alpha_q(q-q_0)] g e_z,
\end{eqnarray}
\begin{eqnarray}
\label{eqa8} \frac{\partial T}{\partial t}+(u\cdot\nabla) T
=Q+\kappa_T \Delta T,
\end{eqnarray}
\begin{eqnarray}
\label{eqa9} \frac{\partial q}{\partial t}+(u \cdot
\nabla)q=G+\kappa_q \Delta q,
\end{eqnarray}
\begin{eqnarray}
\label{eqa10} div u=0,
\end{eqnarray}
where $(\varphi, \theta, z)\in M = S^2_{r_0}\times(r_0, r_0 + h),$
$\nu=\frac{\mu}{\rho_0}$ is the kinematic viscosity, $u=u_\varphi
e_\varphi+u_\theta e_\theta+u_r e_r,$ $(e_\varphi, e_\theta, e_r)$
the local normal basis in the sphereric coordinates, and
$$
\nabla_u u=((u\cdot \nabla) u_\varphi+\frac{u_\varphi
u_z}{r_0}-\frac{u_\varphi u_\theta}{r_0} \tan \theta)e_\varphi+
$$
$$
((u\cdot \nabla) u_\theta+\frac{u_\theta
u_z}{r_0}+\frac{u^2_\varphi}{r_0} \tan \theta)e_\theta+((u\cdot
\nabla) u_z-\frac{u^2_\varphi+ u^2_\theta}{r_0})e_z,
$$
$$
\Delta u=(\Delta u_\varphi+\frac{2}{r_0^2 \cos\theta} \frac{\partial
u_z }{\partial \varphi}+\frac{2 \sin \theta}{r_0^2 \cos^2 \theta}
\frac{\partial u_\theta}{\partial \varphi}- \frac{u_\varphi}{r_0^2
\cos^2 \theta})e_\varphi+
$$
$$
(\Delta u_\theta+\frac{2}{r_0^2} \frac{\partial u_z}{\partial
\theta}-\frac{u_\theta}{r_0^2 \cos^2 \theta}- \frac{2
\sin\theta}{r_0^2 \cos^2 \theta}\frac{u_\varphi}{\partial
\varphi})e_\theta+
$$
$$
 (\Delta u_z+\frac{2 u_0}{r^2}+\frac{2}{r_0^2 \cos \theta}
\frac{\partial (u_\theta \cos\theta)}{\partial \theta}-
\frac{2}{r_0^2 \cos \theta}\frac{u_\varphi}{\partial \varphi})e_z,
$$
$$
\nabla p=\frac{1}{r_0 \cos \theta}\frac{\partial p}{\partial
\theta}e_\varphi+\frac{1}{r} \frac{\partial p}{\partial
\theta}e_\theta+\frac{\partial p}{\partial z}e_z,
$$
$$
div u=\frac{1}{r_0 \cos\theta} \frac{\partial u_\varphi}{\partial
\varphi}+\frac{1}{r_0 \cos\theta} \frac{\partial (u_\theta \cos
\theta)}{\partial \theta}+\frac{\partial u_z}{\partial z},
$$
and the differential operators $(u\cdot \nabla)$ and $\triangle$ are
expressed as
$$
(u \cdot \nabla)=\frac{u_\varphi}{r_0 \cos\theta}
\frac{\partial}{\partial
\varphi}+\frac{u_\theta}{r_0}\frac{\partial}{\partial \theta}+u_z
\frac{\partial}{\partial z},
$$
$$
\Delta=\frac{1}{r_0^2 \cos\theta} \frac{\partial^2}{\partial
\varphi^2}+\frac{1}{r_0^2 \cos\theta} \frac{\partial}{\partial
\theta}(\cos\theta \frac{\partial}{\partial
\theta})+\frac{\partial^2}{\partial z^2}.
$$
\par
Equations (\ref{eqa7})-(\ref{eqa10}) are supplemented with boundary
conditions (\ref{eqa101}), (\ref{eqa104}).

\subsection{Simplification of Model}
Atmospheric circulation is the large-scale motion of the air, which
is essentially a thermal convection process caused by the
temperature and humidity difference between the earth¡¯s surface and
the tropopause. It is a crucial means by which heat and humidity are
distributed on the surface of the earth. Air circulates within the
troposphere, limit vertically by the tropopause at about 8-10km.
Atmospheric motion in the troposphere, together with the oceanic
circulation, plays a crucial role in leading to the global climate
changes and evolution on the earth. There are two types of
circulation cells: the latitudinal circulation and the longitudinal
circulation. The latitudinal circulation is characterized by the
Polar cell, the Ferrel cell, and the Hadley cell, which are major
players in global heat and humidity transport, and do not act alone.
The zonal circulation consists of six circulation cells over the
equator. The overall atmospheric motion is known as the zonal
overturning circulation, and also called the Walker circulation. The
most remarkable feature of the global atmospheric circulation is
that the equatorial Walker circulation divides the whole earth into
three invariant regions of atmospheric flow: the northern
hemisphere, the southern hemisphere, and the equatorial zone. We
also note the important fact that the large-scale structure of the
zonal overturning circulation varies from year to year, but its
basic structure remains fairly stability, it never vanishes. Based
on these natural phenomena, we here present the Zone Hypotheses for
atmospheric dynamics, which amounts to saying that the global
atmospheric system can be divided into three sub-systems: the
North-Hemispheric System, the South-Hemispheric System, and the
Tropical Zone System, which are relatively independent, and have
less influence on each other in their basic structure. More
precisely, the Atmospheric Zone Hypothesis is stated in the
following form\cite{Ma2}.
\par
{\bf  Atmospheric Zone Hypothesis.}  The atmospheric circulation has
three invariant regions: the northern hemisphere domain ($0
<\theta\leq\frac{\pi}{2}$), the southern hemisphere domain ($-
\frac{\pi}{2}\leq \theta < 0$),  and the equatorial zone ($\theta =
0$). Namely, the large-scale circulations in their invariant regions
can act alone with less influence on the others. In particular, the
velocity field $u = (u_\varphi, u_\theta, u_z)$ of the atmospheric
circulation has a vanished latitudinal component in a narrow
equatorial zone, i.e., $u_\theta = 0,$ for
$-\varepsilon<\theta<\varepsilon$, where $\varepsilon> 0$ is a small
number.
\par
 Atmospheric Zone Hypothesis is based on the following
several evidences from theory and practice:
\par
(1) The global atmospheric motion equations
(\ref{eqa7})-(\ref{eqa10}) are of $\theta-$reflexive symmetry, i.e.
under the $\theta-$reflexive transformation $(\varphi, \theta,
z)\rightarrow (\varphi,-\theta, z)$, the velocity field $u$ becomes
$(u_\varphi, u_\theta, u_z)\rightarrow (u_\varphi,-u_\theta, u_z)$,
and equations (\ref{eqa7})-(\ref{eqa10}) are invariant, which
implies that this system is compatible with the Atmospheric Zone
Hypothesis.
\par
(2) Climatic observation data show that when the El
Ni$\tilde{n}o$-Southern Oscillation (the behavior that the Walker
circulation cell in the Western Pacific stops or reverses its
direction) takes place, no oscillation occurs in the latitudinal
cells. It demonstrates the relative independence of these
circulations in their invariant domain.
\par
(3) When a cold current moves southward from Siberia, or a violent
typhoon sweeps northward from the tropics, the weather in Southern
Hemisphere has no response. Atmospheric Zone Hypothesis provides a
theoretic basis for the study of atmospheric dynamics, by which we
can establish locally simplified models to treat many difficult
problems in atmospheric science.
\par
(4) For the three-dimensional atmospheric circulation equation, it
is too difficult to study. So we study the equatorial atmospheric
circulation.
\par
The atmospheric motion equations over the tropics are the equations
restricted on $\theta= 0$, where the meridional velocity component
$u_\theta$ is set to zero, and the effect of the turbulent friction
is taking into considering
\begin{eqnarray}
\label{eqa13}
 \frac{\partial u_\varphi}{\partial t}=-(u \cdot\nabla
) u_\varphi-\frac{u_\varphi u_z}{r_0}+\nu (\Delta u_\varphi+\frac{2
}{r_0^2}\frac{\partial u_z}{\partial
\varphi}-\frac{2u_\varphi}{r_0^2})-\sigma_0 u_\varphi-2 \Omega
u_z-\frac{1}{\rho_0 r_0} \frac{\partial p}{\partial \varphi},
\end{eqnarray}
$$
 \frac{\partial u_z}{\partial t}=-( u\cdot \nabla)
u_z+\frac{u_\varphi^2}{r_0}+\nu (\Delta u_z+\frac{2
}{r_0^2}\frac{\partial u_\varphi}{\partial
\varphi}-\frac{2u_z}{r_0^2})-\sigma_1 u_z-2 \Omega u_\varphi
$$
\begin{eqnarray}
\label{eqa14} -\frac{1}{\rho_0} \frac{\partial p}{\partial
z}-[1-\alpha_T(T-T_0)+\alpha_q(q-q_0)] g,
\end{eqnarray}
\begin{eqnarray}
\label{eqa15} \frac{\partial T}{\partial t}=-(u\cdot\nabla) T
+\kappa_T \Delta T+Q,
\end{eqnarray}
\begin{eqnarray}
\label{eqa16} \frac{\partial q}{\partial t}=-(u \cdot
\nabla)q+\kappa_q \Delta q+G,
\end{eqnarray}
\begin{eqnarray}
\label{eqa17}\frac{1}{r_0}\frac{\partial
u_\varphi}{\partial\varphi}+\frac{\partial u_z}{\partial z}=0,
\end{eqnarray}
Here $\sigma_i = C_i h^2 (i = 0, 1)$ represent the turbulent
friction, $r_0$ is the radius of the earth, the space domain is
taken as $M = S^1_{r_0} \times (r_0, r_0 + h)$ with $S^1_{r_0}$
being the one-dimensional circle with radius $r_0$, and
$$
(u\cdot\nabla)=\frac{u_\varphi}{r_0}\frac{\partial}{\partial
\varphi}+u_z \frac{\partial}{\partial z}, \quad
\Delta=\frac{1}{r_0^2}\frac{\partial^2}{\partial
\varphi^2}+\frac{\partial^2}{\partial z^2}
$$
 For
simplicity, we denote
$$
(x_1,x_2)=(r_0 \varphi, z), \quad (u_1,u_2)=(u_\varphi,u_z).
$$
\par
The atmospheric motion equations(\ref{eqa13})-(\ref{eqa17}) can be
written as
\begin{eqnarray}
\label{eqa18}
 \frac{\partial u_1}{\partial t}=-(u \cdot\nabla
) u_1-\frac{u_1 u_2}{r_0}+\nu (\Delta u_1+\frac{2
}{r_0}\frac{\partial u_2}{\partial x_1}-\frac{2u_1}{r_0^2})-\sigma_0
u_1-2 \Omega u_2-\frac{1}{\rho_0} \frac{\partial p}{\partial x_1},
\end{eqnarray}
$$
 \frac{\partial u_2}{\partial t}=-(u \cdot \nabla)
u_2+\frac{u_1^2}{r_0}+\nu (\Delta u_2+\frac{2 }{r_0}\frac{\partial
u_1}{\partial x_1}-\frac{2u_2}{r_0^2})-\sigma_1 u_2-2 \Omega u_1
$$
\begin{eqnarray}
\label{eqa19} -\frac{1}{\rho_0} \frac{\partial p}{\partial
x_2}-[1-\alpha_T(T-T_0)+\alpha_q(q-q_0)] g,
\end{eqnarray}
\begin{eqnarray}
\label{eqa20} \frac{\partial T}{\partial t}=-(u\cdot\nabla) T
+\kappa_T \Delta T+Q,
\end{eqnarray}
\begin{eqnarray}
\label{eqa21} \frac{\partial q}{\partial t}=-(u \cdot
\nabla)q+\kappa_q \Delta q+G,
\end{eqnarray}
\begin{eqnarray}
\label{eqa22}\frac{\partial u_1}{\partial x_1}+\frac{\partial
u_2}{\partial x_2}=0.
\end{eqnarray}
\par
To make the nondimensional form, let
$$
x=hx^{'}, \quad t=\frac{h^2}{\kappa_T} t^{'} \quad
u=\frac{\kappa_T}{h} u^{'},
$$
$$
T=T_0-(T_0-T_1)\frac{x_2}{h}+(T_0-T_1)T^{'},
$$
$$
q=q_0-(q_0-q_1)\frac{x_2}{h}+(q_0-q_1)q^{'},
$$
$$
p=\frac{\rho_0 \nu \kappa_T
p^{'}}{h^2}-g\rho_0[x_2+\frac{\alpha_T}{2}(T_0-T_1)\frac{x^2_2}{h}-\frac{\alpha_q}{2}(q_0-q_1)\frac{x^2_2}{h}],
$$
\par
The nondimensional form of (\ref{eqa18})-(\ref{eqa22}) reads
$$
 \frac{\partial u^{'}_1}{\partial t^{'}}=-(u^{'} \cdot \nabla
) u^{'}_1+\frac{\nu}{\kappa_T} \Delta u^{'}_1-\frac{\sigma_0
h^2}{\kappa_T} u^{'}_1-\frac{2 \Omega h^2}{\kappa_T}
u^{'}_2-\frac{\nu}{k_T} \frac{\partial p^{'}}{\partial
x^{'}_1}
$$
\begin{eqnarray}
\label{eqa23} -\frac{h}{r_0}u_1u_2+\frac{2h}{r_0 k_T}\frac{\partial
u^{'}_2}{\partial x_1^{'}}-\frac{2h^2}{r_0 \kappa_T} u_1^{'},
\end{eqnarray}
$$
 \frac{\partial u^{'}_2}{\partial t^{'}}=-(u^{'} \cdot \nabla)
u^{'}_2+\frac{\nu}{\kappa_T} \Delta u^{'}_2-\frac{\sigma_1
h^2}{\kappa_T} u^{'}_2-\frac{2 \Omega h^2}{\kappa_T}
u^{'}_1-\frac{h}{r_0}u_1^2+\frac{h}{r_0 k_T}\frac{\partial
u^{'}_1}{\partial x_2^{'}}-\frac{2h^2}{r_0^2 \kappa_T} u_2^{'}
$$
\begin{eqnarray}
\label{eqa24} -\frac{\nu}{\kappa_T} \frac{\partial p^{'}}{\partial
x^{'}_2}+\frac{gh^3}{\kappa_T^2}[\alpha_T(T^0-T_1)T^{'}-\alpha_q(q^0-q_1)q^{'}],
\end{eqnarray}
\begin{eqnarray}
\label{eqa25} \frac{\partial T^{'}}{\partial
t^{'}}=-(u^{'}\cdot\nabla) T^{'} + \Delta
T^{'}+u_2^{'}+\frac{h^2}{(T_0-T_1)\kappa_T}Q,
\end{eqnarray}
\begin{eqnarray}
\label{eqa26} \frac{\partial q^{'}}{\partial t^{'}}=-(u^{'} \cdot
\nabla)q^{'}+\frac{\kappa_q}{\kappa_T} \Delta
q^{'}+u_2^{'}+\frac{h^2}{(T_0-T_1)\kappa_T}G,
\end{eqnarray}
\begin{eqnarray}
\label{eqa27}div u^{'}=0.
\end{eqnarray}
\par
 Because $r_0$ is far larger than
$u_1$, $u_2$,  the atmospheric motion equations
(\ref{eqa23})-(\ref{eqa27}) can be read
\begin{eqnarray}
\label{eqa123}   \frac{\partial u^{'}_1}{\partial t^{'}}=-(u^{'}
\cdot \nabla ) u^{'}_1+\frac{\nu}{\kappa_T} \Delta
u^{'}_1-\frac{\sigma_0 h^2}{\kappa_T} u^{'}_1-\frac{2 \Omega
h^2}{\kappa_T} u^{'}_2-\frac{\nu}{k_T} \frac{\partial
p^{'}}{\partial x^{'}_1},
\end{eqnarray}
$$
 \frac{\partial u^{'}_2}{\partial t^{'}}=-(u^{'} \cdot \nabla)
u^{'}_2+\frac{\nu}{\kappa_T} \Delta u^{'}_2-\frac{\sigma_1
h^2}{\kappa_T} u^{'}_2-\frac{2 \Omega h^2}{\kappa_T} u^{'}_1
$$
\begin{eqnarray}
\label{eqa124} -\frac{\nu}{\kappa_T} \frac{\partial p^{'}}{\partial
x^{'}_2}+\frac{gh^3}{\kappa_T^2}[\alpha_T(T_0-T_1)T^{'}-\alpha_q(q_0-q_1)q^{'}],
\end{eqnarray}
\begin{eqnarray}
\label{eqa125} \frac{\partial T^{'}}{\partial
t^{'}}=-(u^{'}\cdot\nabla) T^{'} + \Delta
T^{'}+u_2^{'}+\frac{h^2}{(T_0-T_1)\kappa_T}Q,
\end{eqnarray}
\begin{eqnarray}
\label{eqa126} \frac{\partial q^{'}}{\partial t^{'}}=-(u^{'} \cdot
\nabla)q^{'}+\frac{\kappa_q}{\kappa_T} \Delta
q^{'}+u_2^{'}+\frac{h^2}{(T_0-T_1)\kappa_T}G,
\end{eqnarray}
\begin{eqnarray}
\label{eqa127}div u^{'}=0.
\end{eqnarray}
\par
Let $P_r=\frac{\nu}{\kappa_T}$, $L_e=\frac{\kappa_q}{\kappa_T}$,
$R=\frac{g\alpha_T(T^0-T_1)h^3}{\kappa_T \nu}$,
$\tilde{R}=\frac{g\alpha_q(q^0-q_1)h^3}{\kappa_T \nu}$,
$\sigma_i^{'}=\frac{\sigma_i h^2}{\nu}$, $\omega=\frac{2 \Omega
h^2}{\nu}$, $Q^{'}=\frac{h^2}{(T_0-T_1)\kappa_T}Q$,
$G^{'}=\frac{h^2}{(T_0-T_1)\kappa_T}G$, omitting the primes, the
nondimensional form of (\ref{eqa123})-(\ref{eqa127}) reads
\begin{eqnarray}
\label{eqa28}
 \frac{\partial u}{\partial t}=P_r(\Delta u-\nabla p-\sigma u)+P_r(RT-\tilde{R}q)\vec{\kappa}-(\nabla \cdot
u) u,
\end{eqnarray}
\begin{eqnarray}
\label{eqa29} \frac{\partial T}{\partial t}=\Delta
T+u_2-(u\cdot\nabla) T + Q,
\end{eqnarray}
\begin{eqnarray}
\label{eqa30} \frac{\partial q}{\partial t}=L_e \Delta q+u_2-(u
\cdot \nabla)q+G,
\end{eqnarray}
\begin{eqnarray}
\label{eqa31}div u=0,
\end{eqnarray}
where $\sigma$ is constant matrix
$$
 \sigma=\left(
   \begin{array}{ll}
\sigma_0 & \omega\\
\omega & \sigma_1
\end{array}
\right).
$$
\par
The problem  (\ref{eqa28})-(\ref{eqa31}) are supplemented with the
following Dirichlet boundary condition at $x_2=0,1$ and periodic
condition for $x_1$:
\begin{eqnarray}
\label{eqa32} (u, T,q ) =0, \quad x_2=0,1,
\end{eqnarray}
\begin{eqnarray}
\label{eqa33} (u, T, q )(0, x_2) = (u, T, q )(2  \pi, x_2),
\end{eqnarray}
and initial value conditions
\begin{eqnarray}
\label{eqa34} (u, T, q ) = (u_0, T_0, q_0), \ \ \ t=0.
\end{eqnarray}

\section{Global Solution of Atmospheric Circulation Equations with \\Humidity Effect}
\setcounter{equation}{0}
\subsection{Preliminaries}

Let $X$ be a linear space, $X_1$, $X_2$ two separable reflexive
Banach spaces, $H$ a Hilbert space. $X_1$, $X_2$ and  $H$ are
completion space of $ X $ under the respective norms. $X_1$,
$X_2\subset H$ are dense embedding. $F:X_2 \times
(0,\infty)\rightarrow X_1^*$ is a continuous mapping. We consider
the abstract equation
\begin{eqnarray}
\label{eqc101} \left\{
   \begin{array}{ll}
\frac{du}{dt}=Fu, &
0<t<\infty, \\
 \\
 u(0)=\varphi,
& \\
\end{array}
\right.
\end{eqnarray}
where $\varphi \in H$, $u: [0,+\infty)\rightarrow H$ is unknown.
\par
{\bf Definition 3.1} Let $\varphi \in H$ be a given initial value.
$u \in L^p((0,T),X_2)\cap L^\infty((0,T), H)$, ($0<T<\infty$) is
called a global solution of Eq(\ref{eqc101}), if $u$ satisfies
$$
(u(t),v)_H=\int_0^t<Fu,v>dt+(\varphi,v)_H, \quad \forall v \in X_1
\subset H.
$$
\par
{\bf Definition 3.2}  Let $u_n, u_0 \in L^p((0,T), X_2)$. $u_n
\rightarrow u_0$ is called uniformly weak convergence in $L^p((0,T),
X_2)$, if $\{u_n\} \subset L^\infty((0,T),H)$ is bounded, and
\begin{eqnarray}
\label{eqc102} \left\{
   \begin{array}{ll}
u_n \rightharpoonup u_0, &\hbox{in} \ \ \ L^p((0,T), X_2), \\
 \\
\lim_{n\rightarrow \infty}\int_0^T |<u_n-u_0,v>_H|^2 dt=0,
& \forall v \in H.\\
\end{array}
\right.
\end{eqnarray}

\par
{\bf Definition 3.3} A mapping $F: X_2 \times (0,\infty)\rightarrow
X_1^*$ is called $T$-weakly continuous, if for $p=(p_1, \ldots,
p_m)$, $0<T<\infty$, and $u_n$ uniformly weakly converge to $u_0$
under Eq(\ref{eqc102}), we have
$$
\lim_{n\rightarrow \infty}\int_0^T <Fu_n,v>dt=\int_0^T <Fu_0,v>dt,
\quad \forall v\in X_1.
$$
\par
{\bf Lemma 3.4} Let $\{u_n\}\subset L^p((0,T), W^{m,p})(m\geq 1)$ be
bounded sequence, and $\{u_n\}$ uniformly weakly converge to $u_0
\subset L^p((0,T), W^{m,p})$, i.e.
\begin{eqnarray}
\label{eqc100} \left\{
   \begin{array}{ll}
u_n \rightharpoonup u_0 \ \ \ \hbox{in}\ \ L^p((0,T), W^{m,p}),& p\geq 2, \\
 \\
\lim_{n\rightarrow \infty}\int_0^T [\int_\Omega (u_n-u_0)v dx
]^2dt=0,
& \forall v \in C_0^\infty(\Omega).\\
\end{array}
\right.
\end{eqnarray}
Then for all $|\alpha|\leq m-1$, we have
$$
D^\alpha u_n \rightarrow D^\alpha u_0\ \ \ \ \hbox{in}\ \ \ \ \
L^2((0,T)\times \Omega).
$$
\par
{\bf Lemma 3.5} $^{\cite{Ma3}}$ Assume $F: X_2 \times
(0,\infty)\rightarrow X_1^*$ is $T$-weakly continuous, and satisfies
\par
(A1) there exists $p=(p_1, \cdots, p_m)$, $p_i>1(1\leq i \leq m)$,
such that
$$
<Fu,u> \leq -C_1\|u\|^p_{X_2}+C_2 \|u\|_H^2+f(t),\quad \forall u\in
X,
$$
where $C_1,C_2$ are constants, $f\in L^1(0,T)(0<T<\infty)$,
$\|\cdot\|^p_{X_2}=\sum_{i=1}^m|\cdot|_i^{p_i}$, $|\cdot|_i$ is
seminorm of $X_2$, $\|\cdot\|_{X_2}=\sum_{i=1}^m |\cdot|_i$,
\par
(A2) there exists $0<\alpha<1$, for all $0<h<1$ and $u \in
C^1([0,\infty),X)$,
$$
|\int_t^{t+h}<Fu,v>dt|\leq Ch^\alpha, \quad \forall v\in X \quad
\hbox{and} \quad 0\leq t<T,
$$
where $C>0$ depends only on $T$, $\|v\|_{X_1}$, $\int_0^t
\|u\|^p_{X_2}dt$ and $\sup_{0\leq t \leq T} \|u\|_H$.\\
 Then for all $\varphi \in H$, Eq(\ref{eqc101}) has a global weak solution
 $$
u\in L^\infty((0,T),H) \cap L^p((0,T), X_2), \quad 0<T<\infty, \quad
p\hbox{\ \ in \ \ (A1)}.
 $$
 \par
 {\bf Remark 3.6} $\|\cdot\|_X$ denotes norm of $X$, and $C_i$ are
 variable constants.

\subsection{Existence of Global Solution}
 We introduce spatial sequences
$$
X=\{\phi=(u,T,q) \in C^\infty(\Omega, R^4)| (u,T,q)  \ \
\hbox{satisfy} \ \ (\ref{eqa31})-(\ref{eqa33})\},
$$
$$
H=\{\phi=(u,T,q) \in L^2(\Omega, R^4)| (u,T,q)\ \   \hbox{satisfy}\
\  (\ref{eqa31})-(\ref{eqa33})\},
$$
$$
H_1=\{\phi=(u,T,q) \in H^1(\Omega, R^4)| (u,T,q) \ \ \hbox{satisfy}\
\  (\ref{eqa31})-(\ref{eqa33})\},
$$

{\bf Theorem 3.7} If $\phi_0=(u_0,T_0,q_0)\in H$, $Q,G \in
L^2(\Omega)$, then Eq(\ref{eqa28})-(\ref{eqa34}) there exists
 a global solution
$$
(u,T,q) \in L^\infty((0,T),H)\cap L^2((0,T),H_1), \quad 0<T<\infty.
$$
\par
{\bf Proof.} Definite $F: H_1 \rightarrow H_1^*$ as
$$
 \begin{array}{lrl}
<F\phi, \psi>&=&\int_\Omega [-P_r \nabla u \nabla v-P_r \sigma u
\cdot v+P_r (RT-\tilde{R}q)v_2-( u\cdot \nabla) u \cdot v
\\\\
&&-\nabla T \nabla S+u_2 S-(u\cdot \nabla)T S+Q S-L_e \nabla q
\nabla z+u_2 z
\\\\
&&-(u\cdot \nabla)q z+G z]dx,\ \ \ \ \hbox{¶Ô}\quad  \forall
\psi=(v,S,z) \in H_1.
 \end{array}
 $$
 \par
 Let $\psi=\phi$. Then
$$
 \begin{array}{lrl}
<F\phi, \phi>&=&\int_\Omega [-P_r |\nabla u|^2-P_r \sigma u \cdot
u+P_r(RT-\tilde{R}q)u_2-(u\cdot\nabla  ) u \cdot u
\\\\
&&-|\nabla T|^2+u_2 T-(u\cdot \nabla)T T +Q T-L_e |\nabla q|^2
\\\\
&&+u_2 q-(u\cdot \nabla)q q+G q]dx
\\\\
&=&\int_\Omega [-P_r |\nabla u|^2-P_r \sigma u \cdot u+(P_r R+1)u_2
T-(P_r \tilde{R}-1)qu_2
\\\\
&&-|\nabla T|^2+ Q T-L_e |\nabla q|^2+G q]dx
\\\\
&\leq&- C_1 \int_\Omega [|\nabla u|^2+|\nabla T|^2+ |\nabla q|^2]dx
+C_2 \int_\Omega[|u|^2+|u||T|+|q||u|
\\\\
&&+|Q| |T|+|G| |q|]dx
\\\\
&\leq& - C_1 \int_\Omega [|\nabla u|^2+|\nabla T|^2+ |\nabla q|^2]dx
+C_2 \int_\Omega[|u|^2+|T|^2+|q|^2]dx
\\\\
&&+C_3 \int_\Omega[|Q|^2+|G|^2]dx
\\\\
&\leq& - C_1 \|\phi\|_{H_1}^2 +C_2 \|\phi\|_{H}^2+C_4,
 \end{array}
 $$
which implies $(A_1)$.
\par
For $\forall \phi \in L^2(0,T), H_1)\cap L^\infty((0,T),H)$ and
$\psi \in X$, $h(0<h<1)$, we have
$$
 \begin{array}{lrl}
&&|\int_t^{t+h}<F\phi, \psi>dt|
\\\\
&=&|\int_t^{t+h} \int_\Omega [-P_r \nabla u \nabla v-P_r \sigma u
\cdot v+P_r (RT-\tilde{R}q)v_2-( u\cdot \nabla) u \cdot v-\nabla T
\nabla S
\\\\
&&+u_2 S-(u\cdot \nabla)T S+ Q S-L_e \nabla q \nabla z+u_2 z-(u\cdot
\nabla)q z+G z]dxdt|
 \end{array}
$$
$$
\begin{array}{lrl}
&\leq& C\int_t^{t+h} \int_\Omega [|\nabla u| |\nabla v|+|u||v|+
|T||v|+|q||v|+|\nabla T| |\nabla S|+|u| |S|
\\\\
&&+|Q| |S| +|\nabla q| |\nabla z|+|u| |z|+|G| |z|]dxdt+C\int_t^{t+h}
[|\int_\Omega ( u\cdot \nabla) u \cdot vdx|
\\\\
&&+|\int_\Omega (u\cdot \nabla)T Sdx|+|\int_\Omega (u\cdot \nabla)q
zdx|]dt
\\\\
&\leq& C\int_t^{t+h}[ (\int_\Omega |\nabla u|^2dx)^{\frac{1}{2}}
(\int_\Omega|\nabla v|^2dx)^{\frac{1}{2}}+
(\int_\Omega|T|^2dx)^{\frac{1}{2}}(\int_\Omega|v|^2dx)^{\frac{1}{2}}
\\\\
&&+(\int_\Omega|q|^2dx)^{\frac{1}{2}}
(\int_\Omega|v|^2dx)^{\frac{1}{2}}+(\int_\Omega|\nabla
T|^2dx)^{\frac{1}{2}} (\int_\Omega|\nabla S|^2dx)^{\frac{1}{2}}
\\\\
&&+(\int_\Omega|u|^2dx)^{\frac{1}{2}}
(\int_\Omega|S|^2dx)^{\frac{1}{2}}+(\int_\Omega|Q|^2dx)^{\frac{1}{2}}
(\int_\Omega|S|^2dx)^{\frac{1}{2}}
\\\\
&&+(\int_\Omega|\nabla q|^2dx)^{\frac{1}{2}} (\int_\Omega|\nabla
z|^2dx)^{\frac{1}{2}}+(\int_\Omega|u|^2dx)^{\frac{1}{2}}
(\int_\Omega|z|^2dx)^{\frac{1}{2}}
\\\\
&&+(\int_\Omega|G|^2dx)^{\frac{1}{2}}
(\int_\Omega|z|^2dx)^{\frac{1}{2}}]dt+C\int_t^{t+h}
[|\sum_{i,j=1}^2\int_\Omega ( u_i u_j \frac{\partial v_j}{\partial
x_i}dx|
\\\\
&&+|\sum_{i=1}^2\int_\Omega u_iT \frac{\partial S}{\partial
x_i}dx|+|\sum_{i=1}^2\int_\Omega u_iq\frac{\partial z}{\partial
x_i}dx|]dt
\\\\
&\leq& C(\|u\|_{L^2(0,T;H^1)} \|D v\|_{L^2}h^{\frac{1}{2}}+
\|T\|_{L^2(0,T;L^2)}\|v\|_{L^2}h^{\frac{1}{2}}+\|q\|_{L^2(0,T;L^2)}
\|v\|_{L^2}h^{\frac{1}{2}}
\\\\
&&+\|T\|_{L^2(0,T;H^1)} \|D
S\|_{L^2}h^{\frac{1}{2}}+\|u\|_{L^\infty(0,T;L^2)}\|S\|_{L^2}h^{\frac{1}{2}}
+\|Q\|_{L^2}\|S\|_{L^2}h
\\\\
 &&+\|
q\|_{L^2(0,T;H^1)} \|D
z\|_{L^2}h^{\frac{1}{2}}+\|u\|_{L^\infty(0,T;H)}
\|z\|_{L^2}h+\|G\|_{L^2}\|z\|_{L^2}h
\\\\
&&+\|v\|_{C^1}\|u\|_{L^\infty(0,T,L^2)}h+\|S\|_{C^1}
\|u\|_{L^\infty(0,T;L^2)}^{\frac{1}{2}}\|T\|_{L^\infty(0,T;L^2)}^{\frac{1}{2}}h
\\\\
&&+\|z\|_{C^1}
\|u\|_{L^\infty(0,T;L^2)}^{\frac{1}{2}}\|q\|_{L^\infty(0,T;L^2)}^{\frac{1}{2}}h)
\\\\
&\leq& C h^{\frac{1}{2}},
 \end{array}
$$
which implies $(A_2)$.
\par
We will prove that $F: H_1 \rightarrow H_1^*$ is $T$-weakly
continuous. Let $\phi_n=(u^n,T^n,q^n) \rightharpoonup \phi_0
=(u^0,T^0,q^0)$ be uniformly convergence, i.e., $\{\phi_n\} \subset
L^\infty((0,T);H)$ is bounded, and
 $$
 \left\{
   \begin{array}{lcl}
\phi_n \rightharpoonup \phi_0 & \hbox{in} & L^p((0,T);H_1),
 \\
\lim_{n\rightarrow \infty}
\int_0^T|<\phi_n-\phi_0,\psi>_H|^2dt=0,&&\forall \psi \in H.
\end{array}
\right.
$$
\par
From Lemma 3.4, we known that $\phi_n \rightarrow \phi_0$ in
$L^2(\Omega \times (0,T))$.
\par
Then $\forall \psi \in X \subset C^\infty(\Omega, R^4)\cap H_1$, we
have
$$
\lim_{n \rightarrow \infty}\int_0^t\int_\Omega(u^n \cdot \nabla)u^n
\cdot v dxdt=\lim_{n \rightarrow
\infty}\int_0^t\int_\Omega\sum_{i,j=1}^2 u^n_i \frac{\partial
u^n_j}{\partial x_i} v_jdxdt
$$
$$
=-\lim_{n \rightarrow \infty}\int_0^t\int_\Omega\sum_{i,j=1}^2 u^n_i
u^n_j\frac{\partial v_j}{\partial x_i} dxdt
$$
$$
=-\int_0^t\int_\Omega(u^0 \cdot \nabla)v  \cdot u^0dxdt
$$
$$
=\int_0^t\int_\Omega(u^0 \cdot \nabla)u^0 \cdot v dxdt,
$$
and
$$
\lim_{n \rightarrow \infty}\int_0^t\int_\Omega(u^n \cdot \nabla)T^n
S dxdt=\lim_{n \rightarrow \infty}\int_0^t\int_\Omega\sum_{i=1}^2
u^n_i \frac{\partial T^n}{\partial x_i} Sdxdt
$$
$$
=-\lim_{n \rightarrow \infty}\int_0^t\int_\Omega\sum_{i=1}^2 u^n_i
T^n\frac{\partial S}{\partial x_i} dxdt
$$
$$
=-\int_0^t\int_\Omega(u^0 \cdot \nabla)S   T^0dxdt
$$
$$
=\int_0^t\int_\Omega(u^0 \cdot \nabla)T^0 S dxdt,
$$
and
$$
\lim_{n \rightarrow \infty}\int_0^t\int_\Omega(u^n \cdot \nabla)q^n
z dxdt=\lim_{n \rightarrow \infty}\int_0^t\int_\Omega\sum_{i=1}^2
u^n_i \frac{\partial q^n}{\partial x_i} zdxdt
$$
$$
=-\lim_{n \rightarrow \infty}\int_0^t\int_\Omega\sum_{i=1}^2 u^n_i
q^n\frac{\partial z}{\partial x_i} dxdt
$$
$$
=-\int_0^t\int_\Omega(u^0 \cdot \nabla)z   q^0dxdt
$$
$$
=\int_0^t\int_\Omega(u^0 \cdot \nabla)q^0 z dxdt,
$$
Thus,
\begin{eqnarray}
\label{eqc9} \lim_{n \rightarrow \infty}\int_0^t <F\phi_n,
\psi>dt=\int_0^t <F\phi_0, \psi>dt.
\end{eqnarray}
\par
Because $X$ is dense in $H_1$, Eq(\ref{eqc9}) holds for $\psi \in
H_1$. In other words, the mapping $F: H_1 \rightarrow H_1^*$ is
$T$-weakly continuous.
\par
From Lemma 3.5, Eq(\ref{eqa28})-(\ref{eqa34}) has a global weak
solution
$$
(u,T,q) \in L^\infty((0,T),H)\cap L^2((0,T),H_1), \quad 0<T<\infty.
$$
$\Box$
\par
{\bf Remark 3.8} Existence of global solutions to the atmospheric
circulation models implies that atmospheric circulation has its own
running way as humidity source and heat source change, and confirms
that the atmospheric circulation models are reasonable.

\end{document}